\documentclass[11pt,a4paper]{amsart}
\usepackage{graphicx,multirow,array,amsmath,amssymb}

\newtheorem{theorem}{Theorem}
\newtheorem{proposition}[theorem]{Proposition}

\begin{document}

\title[Min. lattice length and ropelength of 2-bridge knots]{Minimum lattice length and ropelength of 2-bridge knots and links}
\author[Y. Huh]{Youngsik Huh}
\address{Department of Mathematics, College of Natural Sciences, Hanyang University, Seoul 133-791, Korea}
\email{yshuh@hanyang.ac.kr}
\author[K. Hong]{Kyungpyo Hong}
\address{Department of Mathematics, Korea University, Anam-dong, Sungbuk-ku, Seoul 136-701, Korea}
\email{cguyhbjm@korea.ac.kr}
\author[H. Kim]{Hyoungjun Kim}
\address{Department of Mathematics, Korea University, Anam-dong, Sungbuk-ku, Seoul 136-701, Korea}
\email{kimhjun@korea.ac.kr}
\author[S. No]{Sungjong No}
\address{Department of Statistics, Ewha Womans University, 52, Ewhayeodae-gil, Seodaemun-gu, Seoul 120-750, Korea}
\email{sungjongno84@gmail.com}
\author[S. Oh]{Seungsang Oh}
\address{Department of Mathematics, Korea University, Anam-dong, Sungbuk-ku, Seoul 136-701, Korea}
\email{seungsang@korea.ac.kr}

\thanks{This corresponding author(Seungsang Oh) was supported by Basic Science Research Program through
the National Research Foundation of Korea(NRF) funded by the Ministry of Science,
ICT \& Future Planning(MSIP) (No.~2011-0021795).}

\begin{abstract}
Knots are commonly found in molecular chains such as DNA and proteins,
and they have been considered to be useful models for structural analysis of these molecules.
One interested quantity is the minimum number of monomers
necessary to realize a molecular knot.
The minimum lattice length $\mbox{Len}(K)$ of a knot $K$ indicates
the minimum length necessary to construct $K$ in the cubic lattice.
Another important quantity in physical knot theory is the ropelength
which is one of knot energies measuring the complexity of knot conformation.
The minimum ropelength $\mbox{Rop}(K)$ is the minimum length of an ideally flexible rope
necessary to tie a given knot $K$.

Much effort has been invested in the research project for finding upper bounds on both quantities
in terms of the minimum crossing number $c(K)$ of the knot.
It is known that $\mbox{Len}(K)$ and $\mbox{Rop}(K)$ lie between
$\mbox{O}(c(K)^{\frac{3}{4}})$ and $\mbox{O}(c(K) [\ln (c(K))]^5)$,
but unknown yet whether any family of knots has superlinear growth.
In this paper, we focus on 2-bridge knots and links.
Linear growth upper bounds on the minimum lattice length and minimum ropelength
for nontrivial 2-bridge knots or links are presented:

\vspace{2mm}
\hspace{5mm} $\mbox{Len}(K) \leq 8 c(K) + 2$.

\vspace{2mm}
\hspace{5mm} $\mbox{Rop}(K) \leq 11.39 c(K) + 12.37$.

\end{abstract}

\maketitle

\section{Introduction} \label{sec:intro}

A knot is an embedding of a circle in 3-dimensional Euclidean space,
and a link is a disjoint collection of knots.
Knots have been considered to be useful models for simulating molecular chains such as DNA and proteins.
Especially the length of the polymer is one of key parameters
which impact on the topology of a macromolecule.
In this paper, we consider two kinds of measures of the complexity of knot conformation
related to the length of knots.

We first consider a knot in the cubic lattice
$\mathbb{Z}^3=(\mathbb{R} \times \mathbb{Z} \times \mathbb{Z}) \cup (\mathbb{Z}
\times \mathbb{R} \times \mathbb{Z}) \cup (\mathbb{Z} \times \mathbb{Z} \times \mathbb{R})$.
An {\em edge\/} is a line segment of unit length joining two nearby lattice points in $\mathbb{Z}^3$.
An edge parallel to the $x$-axis is called an $x$-{\em edge\/},
and the plane with the equation $x=i$ for some integer $i$ is called $x$-{\em level\/} $i$.
The terminologies concerning the $y$- and $z$-coordinates are defined in the same manner.
The minimum number of edges necessary to construct a given knot $K$ in $\mathbb{Z}^3$
is called the {\em minimum lattice length\/}, denoted by $\mbox{Len}(K)$.
Diao \cite{D1} introduced this terminology (he used ``minimal edge number" instead) and
proved that the minimal lattice length of the trefoil knot $3_1$ is $24$.
This kind of polygonal representations of knots are very useful for many applications in science.
The microscopic level molecules are more similar to rigid sticks than flexible ropes.
In fact, DNA molecules are made up of small rigid sticks of sugar,
phosphorus, nucleotide proteins and hydrogen bonds.
This quantity is closely related to the minimum number of monomers
necessary to realize a molecular knot.

Another important quantity in physical knot theory is the ropelength
which is one of knot energies measuring the complexity of knot conformation.
Minimum energy conformations are considered canonical or ideal conformations.
The ropelength of a knot is the quotient of its length by its thickness,
the radius of the largest embedded normal tube around the knot.
The minimum ropelength of a knot $K$ is denoted by $\mbox{Rop}(K)$.
The ropelength of a knot was defined in \cite{BO} and the basic theory was developed in \cite{LSDR}.
This shortest tube of uniform thickness forming a given knot represents
the canonical or ideal geometric representation of the knot.
Ideal knots provide irreducible representations of knots
which are related to physical features such as the time-averaged shapes of
knotted DNA molecules in solution.

Both quantities are closely correlated such as $\mbox{Rop}(K) \leq 2 \ \mbox{Len}(K)$.
Much effort has been invested in the research project for finding upper bounds on these quantities
in terms of the minimum crossing number $c(K)$ of the knot.
It is already known that $\mbox{Len}(K)$ and $\mbox{Rop}(K)$ lie between
$\mbox{O}(c(K)^{\frac{3}{4}})$ and $\mbox{O}(c(K) [\ln (c(K))]^5)$,
but unknown yet whether any family of knots has superlinear growth \cite{B,DEPZ}.

In this paper, we focus on 2-bridge knots or links.
Linear growth upper bounds on the minimum lattice length and minimum ropelength for
nontrivial 2-bridge knots or links are presented:

\begin{theorem} \label{thm:Len}
Let $K$ be a nontrivial 2-bridge knot or link.
Then,

\vspace{2mm}
\hspace{5mm} ${\rm Len}(K) \leq 8 c(K) + 2.$
\end{theorem}

\begin{proposition} \label{prop:Rop}
Let $K$ be a nontrivial 2-bridge knot or link.
Then, for any constant $h \geq 1.205$,

\vspace{2mm}
\hspace{5mm} ${\rm Rop}(K) \leq 2h(\sqrt{\pi^2+4}+1) c(K) + 4 \pi + 14h.$
\end{proposition}

By setting $h=1.205$, we have an upper bound of the minimum ropelength:
${\rm Rop}(K) \leq 11.39 c(K) + 29.44.$

In Section 3 and 4, we construct 2-bridge knots or links in specific ways
which realize the upper bounds in Theorem \ref{thm:Len} and Proposition \ref{prop:Rop}, respectively.
In Section 5, the 2-bridge knots or links constructed in Section 4 are locally modified
so that the constant term of the upper bound is improved as follows.

\begin{theorem} \label{thm:Rop2}
Let $K$ be a nontrivial 2-bridge knot or link with $c(K) \geq 6$.
Then,

\vspace{2mm}
\hspace{5mm} ${\rm Rop}(K) \leq 11.39 c(K) + 12.37.$
\end{theorem}

In fact, a linear growth upper bound on the minimum ropelength for
a much larger class of knots called Conway algebraic knots is known \cite{DEZ}.

For some knots with small minimum crossing numbers,
the exact values of the minimum lattice length were mathematically confirmed.
Diao \cite{D1} proved rigorously that the minimum lattice length of the trefoil knot $3_1$ is 24
and all the other nontrivial knots need more than 24 edges.
Scharein {\em et al\/} \cite{SIADSV} proved that
the minimum lattice length of $4_1$ and $5_1$ are 30 and 34, respectively.
The reader can find numerical estimations for various knots in \cite{HNRAV, JaP, SIADSV}.

Many results about finding lower bounds for the ropelength appear in \cite{B, BS, CKS1, CKS2, D2, DE,DET}.
Cantarella {\em et al\/} \cite{CFM} found an upper bound of the minimum ropelength of a knot
or non-split link:

\vspace{2mm}
\hspace{5mm} $\mbox{Rop}(K) \leq 1.64 c(K)^2 + 7.69 c(K) + 6.74.$
\vspace{2mm}

Diao {\em et al\/} \cite{DEPZ, DEY} established $\mbox{O}(c(K)^{\frac{3}{2}})$
(and later $\mbox{O}(c(K) [\ln (c(K))]^5)$)
upper bounds for the minimum lattice length and the minimum ropelength:

\vspace{2mm}
\hspace{5mm} $\mbox{Len}(K) \leq 136 c(K)^{\frac{3}{2}} + 84 c(K) + 22 c(K)^{\frac{1}{2}} + 11.$

\vspace{2mm}
\hspace{5mm} $\mbox{Rop}(K) \leq 272 c(K)^{\frac{3}{2}} + 168 c(K) + 44 c(K)^{\frac{1}{2}} + 22.$

\vspace{2mm}
\hspace{5mm} $\mbox{Len}(K), \ \mbox{Rop}(K) \leq \mbox{O}(c(K) [\ln (c(K))]^5).$
\vspace{2mm}

Recently the authors \cite{HKNO} found other $\mbox{O}(c(K)^2)$ upper bounds with smaller coefficients:

\vspace{2mm}
\hspace{5mm} $\mbox{Len}(K) \leq \min \left\{ \frac{3}{4}c(K)^2 + 5c(K) + \frac{17}{4}, \,
\frac{5}{8}c(K)^2 + \frac{15}{2}c(K) + \frac{71}{8} \right\}.$

\vspace{2mm}
\hspace{5mm} $\mbox{Rop}(K) \leq \min \left\{ \hspace{-1mm}
\begin{tabular}{l}
$1.5 c(K)^2  + 9.15 c(K) + 6.79$, \\
$1.25 c(K)^2  + 14.58 c(K) + 16.90$
\end{tabular}
\hspace{-1mm} \right\}.$

\section{Standard diagrams of 2-bridge knots and links}

\begin{figure}[h]
\includegraphics{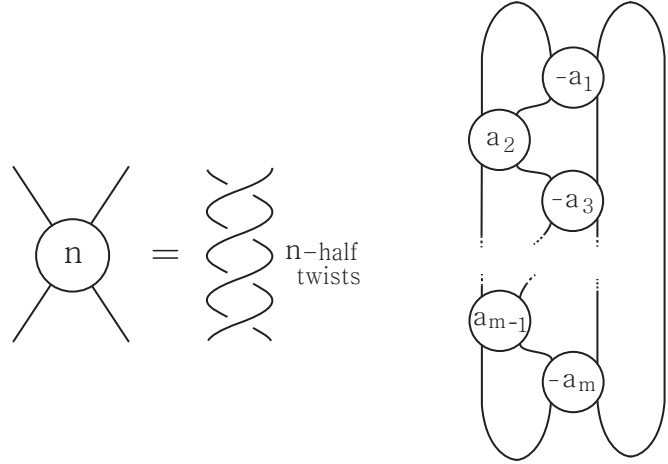}
\caption{Conway notation $(a_1,a_2,\cdots,a_m)$}
\label{fig1}
\end{figure}

In this section we briefly review the standard diagram of a $2$-bridge knot or link
in terms of the Conway notation.
Conway \cite{C} introduced the concept of a tangle in a knot or link diagram
which is a region in the diagram surrounded by a circle
such that the knot or link crosses the sphere exactly four times.
An integral tangle is made from two strands that wrap around each other,
and identified by the number of half-twists within it.
More precisely the integer inside the circle is positive
if it indicates the number of right-handed half-twists and negative if left-handed,
as in Figure \ref{fig1}.
These integral tangles are connected together as the right figure to form
a 2-bridge knot or link which is represented by a Conway notation $(a_1,a_2,\cdots,a_m)$.
Note that if all $a_i$ are positive integers,
then the positive and negative signs of integers in the figure appear alternately,
so it gives a non-nugatory alternating diagram of a 2-bridge knot or link.

As summarized in \cite[Section 2]{Mc},
any nontrivial 2-bridge knot or link can be represented by Conway notation
$(a_1,a_2,\cdots,a_m)$ with positive integers $a_i$ and odd number $m$
due to work by Burde and Zieschang \cite{BZ},
and this non-nugatory alternating diagram gives the minimum number of crossings
due to Kauffman \cite{K}, Murasugi \cite{Mu} and Thistlethwaite \cite{T}.

\section{Minimum lattice length of 2-bridge knots or links}

In this section we prove Theorem \ref{thm:Len}.
Let $K$ be a nontrivial knot or link whose bridge number is 2.
Consider a standard diagram of $K$ in terms of the Conway notation $(a_1,a_2,\cdots,a_m)$
where all $a_i$ are positive and $m$ is an odd number.
Remark that this diagram has explicit $c(K)$ crossings which is $a_1 + \cdots + a_m$.

\vspace{3mm}
\noindent {\bf Step 1.}
{\em Embedding of a 2-bridge knot or link $K$ into the cubic lattice $\mathbb{Z}^3$.\/}
\vspace{2mm}

\begin{figure}[h]
\includegraphics{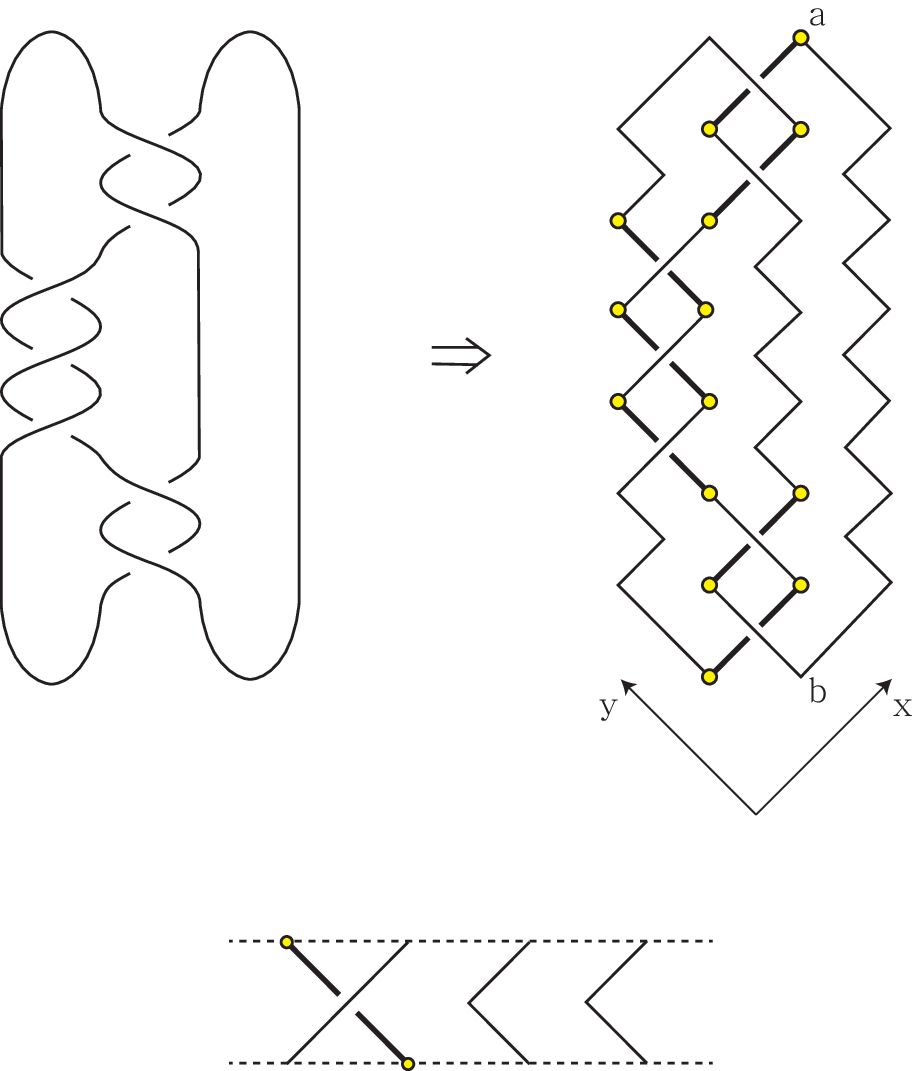}
\caption{2-bridge knot or link in the cubic lattice}
\label{fig2}
\end{figure}

We settle this diagram of $K$ into $\mathbb{Z}^3$.
See Figure \ref{fig2} for an example of a 2-bridge link with the Conway notation $(2,3,2)$ with 7 crossings.
For better view, we rotate the axes of coordinates $45^{\circ}$ counterclockwise.
The $z$-axis is perpendicular to the $xy$-plane.
First draw the diagram on the plane which is on the $z$-level 2
so that it consists of $x$-edges and $y$-edges.
The bold edges are put on the $z$-level 1.
And bold dots represent $z$-edges between the $z$-levels 1 and 2.
Then we have a lattice presentation of $K$ in $\mathbb{Z}^3$.

We can easily count the number of all edges used in this construction.
4 $x$-edges, 4 $y$-edges and 2 $z$-edges are needed for each floor representing one crossing
as drawn in the bottom figure.
At the top and the bottom floors, the same number of edges are needed.
This implies that we need $10 c(K)$ edges in total.

\vspace{3mm}
\noindent {\bf Step 2.}
{\em Folding argument to reduce one fifth of edges.\/}
\vspace{2mm}

\begin{figure}[h]
\includegraphics{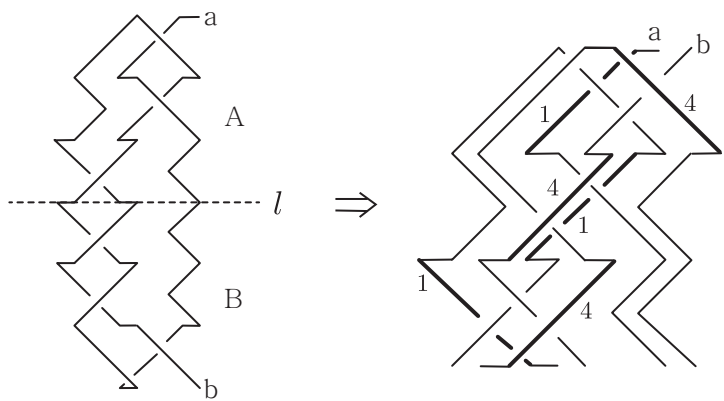}
\caption{Folding argument}
\label{fig3}
\end{figure}

Fold the lattice knot or link diagram in Step 1 to reduce one fifth of the number of edges as follows.
First consider the rightmost arc between the vertices $a$ and $b$ on the $z$-level 2
which consists of $2c(K)$ edges.
Delete the arc and split the remaining part into two pieces, named $A$ and $B$,
by cutting at three vertices on the $z$-level 2 lying on a middle line $l$ as drawn in Figure \ref{fig3}.

This middle line indicates a line $y=x+k$ for some integer $k$ on the $xy$-plane
such that if $c(K)$ is even, then $l$ locates exactly at the center of the height
as viewed of a 2-bridge diagram,
but if $c(K)$ is odd, then $l$ locates at a little above the center of the height passing bold dots.
Note that these three cutting points are not from crossing points of the diagram.
Rotate the bottom piece $B$ by $180^{\circ}$ around $l$ and push up $B$ into the $z$-levels 3 and 4.

\begin{figure}[h]
\includegraphics{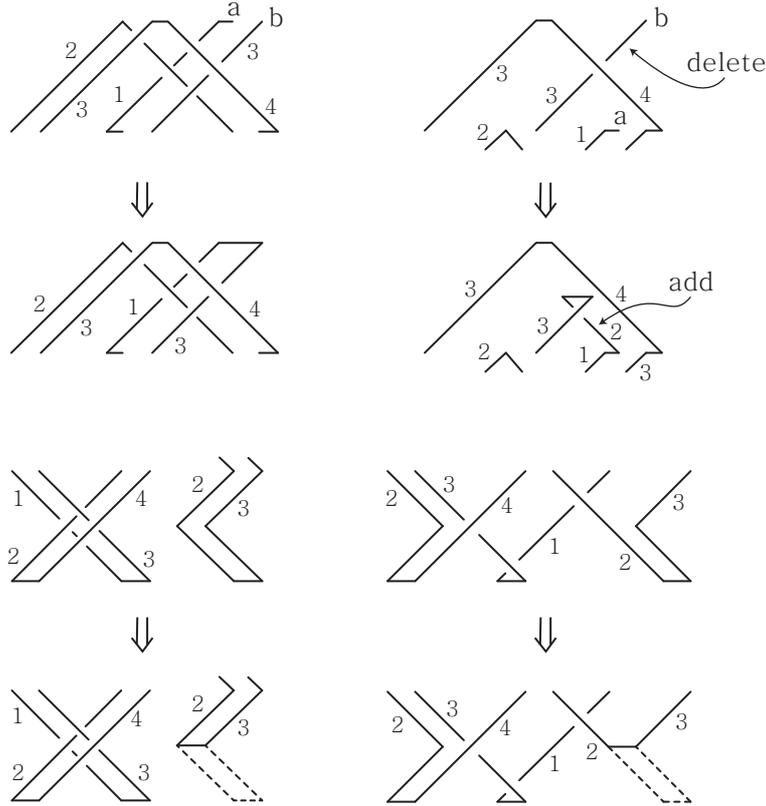}
\caption{Connecting four pairs of cutting points}
\label{fig4}
\end{figure}

If $c(K)$ is even, connect $a$ and $b$ by adding one $z$-edge,
and if $c(K)$ is odd, first delete one $x$-edge adjacent to $b$ and
next add one $y$-edge near $a$ before adding one $z$-edge as shown in Figure \ref{fig4}.
Now connect the other three pairs of the cutting points near $l$ by three $z$-edges.
We still have chance to reduce two more $x$-edges (or $y$-edges) in any case
as illustrated in two bottom figures.

We count the number of edges.
By deleting the rightmost arc, the number of edges are reduced by $2c(K)$.
Then we add 4 $z$-edges to connect four pairs of cutting points
and subtract 2 $x$-edges (or $y$-edges).
This guarantees that $8 c(K) + 2$ edges are enough.
This completes the proof.

\section{Minimum ropelength of 2-bridge knots or links}

In this section we prove Proposition \ref{prop:Rop}.
The basic strategy for the proof is similar to the proof of Theorem \ref{thm:Len}.
Let $K$ be a nontrivial 2-bridge knot or link with a standard diagram associated to the Conway notation $(a_1,a_2,\cdots,a_m)$.
We assume that the rope of $K$ has radius 1 everywhere for simple calculation.

\vspace{3mm}
\noindent {\bf Step 1.}
{\em Embedding of a 2-bridge knot or link $K$ into three cylindrical towers.\/}
\vspace{2mm}

\begin{figure}[h]
\includegraphics{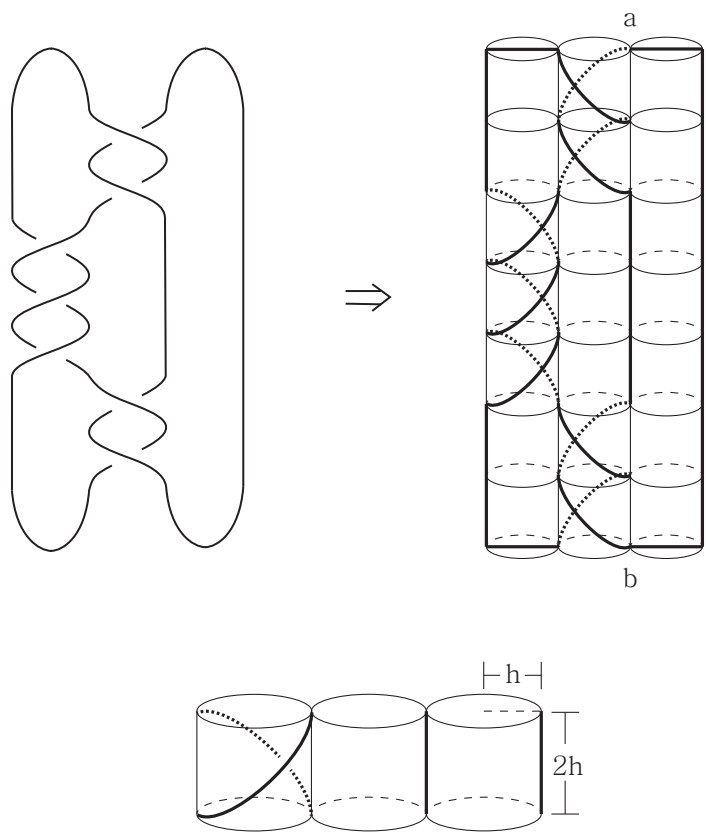}
\caption{2-bridge knot or link in three cylindrical towers}
\label{fig5}
\end{figure}

A cylindrical tower is the stack of cylinders such that the radius and the height
of each cylinder are $h$ and $2h$ for a real number $h \geq 1$.
First draw the diagram of $K$ on three parallel consecutive cylindrical towers as illustrated in Figure \ref{fig5}.
In each floor, exactly one of three cylinders is associated with a crossing of the diagram.

We can easily calculate the length of this embedding.
In each floor, there are two helical arcs with length $h \sqrt{\pi^2 + 4}$
and two vertical line segments with length $2h$ as drawn in the bottom figure.
Additionally, at the top and the bottom,
we need four more horizontal line segments with length $2h$ which come from the diameters of cylinders.
Therefore the length of the embedding is $2h(\sqrt{\pi^2 + 4}+2)c(K) + 8h$.

\vspace{3mm}
\noindent {\bf Step 2.}
{\em Folding argument to reduce about $2hc(K)$.\/}
\vspace{2mm}

\begin{figure}[h]
\includegraphics{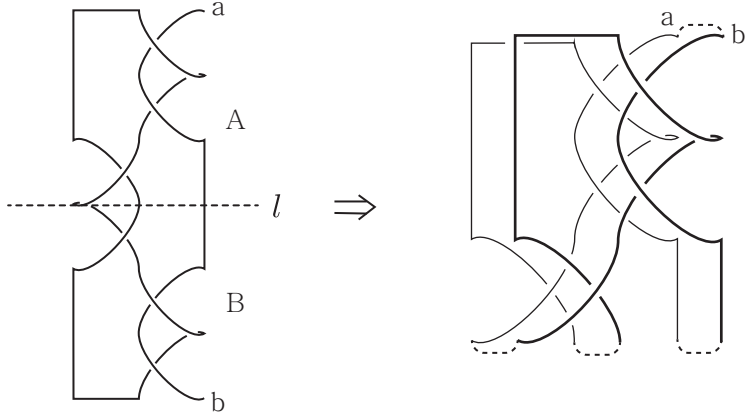}
\caption{Folding argument}
\label{fig6}
\end{figure}

Fold the cylindrical towers to reduce the length of the embedding as follows.
First consider the rightmost arc between vertices $a$ and $b$ which has the length $2h c(K) + 4h$.
Delete the arc and split the remaining part into two pieces, named $A$ and $B$,
by cutting these cylindrical towers at a middle level $l$ of the height as drawn in Figure \ref{fig6}.
If $c(K)$ is even, then $l$ locates exactly at the middle of the height as viewed of a 2-bridge diagram,
but if $c(K)$ is odd, then $l$ locates at the bottom level of three cylinders
which are located at the middle of the height.
Note that these three cutting points are not from crossing points of the diagram.
Rotate the bottom piece $B$ by $180^{\circ}$ around $l$ and pull it in front of the piece $A$.
We make a space of distance 2 between the cylindrical towers $A$
and the cylindrical towers $B$.

\begin{figure}[h]
\includegraphics{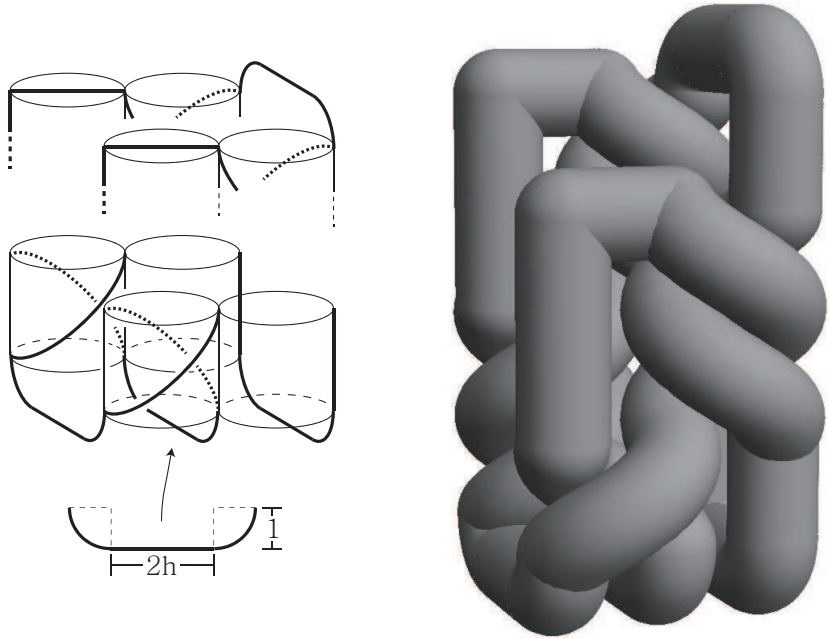}
\caption{Connecting four pairs of cutting points}
\label{fig7}
\end{figure}

Now we connect three pairs of cutting points at level $l$ and another pair of $a$ and $b$
by adding four arcs as drawn in Figure \ref{fig7}.
The arc consists of two quarter circles of radius 1 and a line segment of length $2h$ connecting them.
Note that if $c(K)$ is odd, then we first attach a vertical line segment with length $2h$
to $b$ so that the other endpoint of this segment has the same height as $a$
before adding the connecting arc.
The right figure shows the tube link with uniform radius 1
resulted from the 2-bridge link in Figure \ref{fig5}.

Finally we measure the total length of the embedding of $K$.
We need two helical arcs with length $h \sqrt{\pi^2 + 4}$
and a vertical line segment with length $2h$ for each crossing,
two segments with length $2h$ at the top, and four connecting arcs with length $\pi + 2h$.
If $c(K)$ is odd, then we need one more vertical line segment with length $2h$ as mentioned above.
This implies that the length of this embedding is bounded above by $2h (\sqrt{\pi^2 + 4}+1)c(K) + 4 \pi +14h$.

\vspace{3mm}
\noindent {\bf Step 3.}
{\em Finding a lower bound of $h$.\/}
\vspace{2mm}

\begin{figure}[h]
\includegraphics{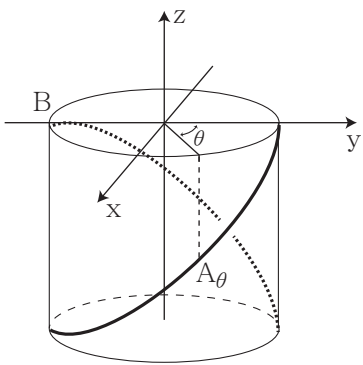}
\caption{Distance between two arcs on a cylinder}
\label{fig8}
\end{figure}

We find a proper lower bound of $h$ avoiding that the rope overlaps itself.
Obviously $h$ should be greater than 1.
Consider the shortest distance between two arcs in a cylinder representing a crossing.
Let $B$ be an end point of the under-crossing arc which intersects the top of the cylinder,
and $A$ be the over-crossing arc as in Figure \ref{fig8}.

Indeed, it is enough to check the distance from the point $B$ to the arc $A$ because of the symmetry.
Put the cylinder into $\mathbb{R}^3$
so that the coordinate of the center of top disk is the origin, and the coordinate of $B$ is $(0,-h,0)$.
So $A$ can be parametrized
as $A_{\theta} = (h \sin \theta, h \cos \theta, -2h \frac{\theta}{\pi})$ for $\theta \in [0 , \pi ]$
where the angle $\theta$ is shown in the figure.
The distance between $B$ and a point $A_{\theta}$ is
$$\mbox{dist} \{B, A_{\theta} \}  = h \sqrt{2 + 2 \cos \theta + 4 \frac{\theta^2}{\pi^2}}.$$
Taylor's theorem gives the following approximation for $\theta \in [0 , \pi ]$:
$$2 + 2 \cos \theta + 4 \frac{\theta^2}{\pi^2} >
4 + (\frac{4}{\pi^2} - 1) \theta^2 + \frac{2}{4!} \theta^4 - \frac{2}{6!} \theta^6.$$
Let $f(\theta)$ be the function on the right side of the inequality.
Then,
$$f'(\theta) = - \frac{2}{5!} \theta (\theta^4 - 20 \theta^2 + 5! (1 - \frac{4}{\pi^2}))$$
and $\theta_{\circ} = \sqrt{10 - \sqrt{100 - 120 (1 - \frac{4}{\pi^2})}}  \approx 2.3946$
is the unique root of the equation $f'(\theta) = 0$ in $(0, \pi)$.
Indeed $f(\theta)$ has the minimum value at $\theta_{\circ}$ among $[0 , \pi ]$.
Since $\mbox{dist} \{B, A_{\theta} \} > h \sqrt{f(\theta)}$, the condition $h \sqrt{f(\theta_{\circ})} \geq 2$
guarantees that $\mbox{dist} \{B, A_{\theta} \}$ is greater than 2 for all $\theta \in [0 , \pi ]$.
To satisfy this condition,
$$h \geq \frac{2}{\sqrt{f(\theta_{\circ})}} \approx 1.2045.$$

This completes the proof of Proposition \ref{prop:Rop}.

\section{Reduction of the constant term}

In this section we prove Theorem \ref{thm:Rop2}.
Let $K$ be a nontrivial 2-bridge knot or link constructed in
a folded cylindrical towers through the procedure in Section 4.
We assume that $c(K) \geq 6$.
We modify $K$ at the bottom and the top parts to reduce the constant term of
the upper bound of the minimum ropelength $11.39 c(K) + 29.44$
obtained from Proposition \ref{prop:Rop} by setting $h=1.205$.

Divide these cylindrical towers and $K$ into three parts,
named the bottom, the middle and the top parts.
The bottom (and the top) part indicates four cylinders at the bottom (and the top)
and three arcs of $K$ lying on or below (and above, respectively) these cylinders.
The middle part indicates the rest $2  (c(K)-4)$ cylinders between them
and six subarcs of $K$ lying on the cylinders.
Recall that the ropelengths of the bottom, top and middle parts obtained in Step 2 of Section 4 are
$4h (\sqrt{\pi^2 + 4}+1) + 3 \pi +6h$ ($\approx 39.43$),
at most $4h (\sqrt{\pi^2 + 4}+1) + \pi +8h$ ($\approx 35.55$),
and $2h (\sqrt{\pi^2 + 4}+1)(c(K)-4)$ ($\approx 11.39 c(K) - 45.54$), respectively.
\vspace{3mm}

\noindent {\bf Step 1.}
{\em Shortening the bottom part.\/}
\vspace{2mm}

\begin{figure}[h]
\includegraphics{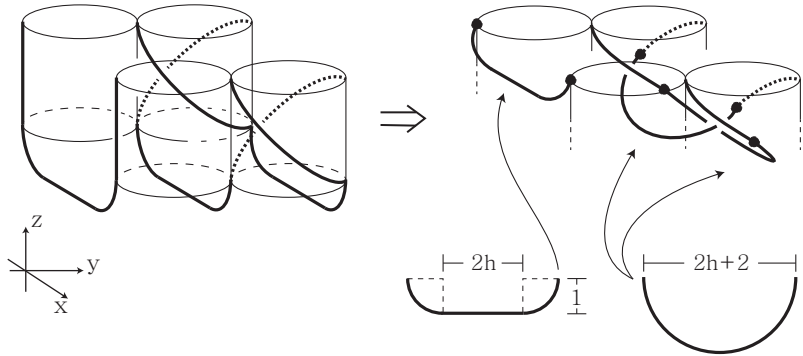}
\caption{Two crossings lying on the same side at the bottom part}
\label{fig9}
\end{figure}

For the bottom part, there are four possible types according to the positions of two crossings.
First consider the type that the two crossings lie on the two right cylinders as shown in Figure \ref{fig9}.
We have indeed the same result for the type that they lie on the two left cylinders
because two types are merely mirror reflections of each other with respect to an $xz$-plane.
Replace the bottom part by the three arcs illustrated in the right figure.
In this case, the only four upper half helical arcs associated to the crossings are kept,
and the other parts are deleted.
Instead, attach two half circles of radius $h+1$ to two pairs of end-points of the half helical arcs.
To maintain the distance between any pair of the attached half circles at least 2,
the tangent line at an end-point of each half circles forms the angle $45^{\circ}$ downward
with respect to the $z$-axis.
Also attach an arc consisting of two quarter circles of radius 1 and a horizontal line segment
of length $2h$ connecting them.
The total length of three newly constructed arcs is
$2h\sqrt{\pi^2+4} + 2 \pi (h+1) + (\pi +2h)$ ($\approx 28.38$).

\begin{figure}[h]
\includegraphics{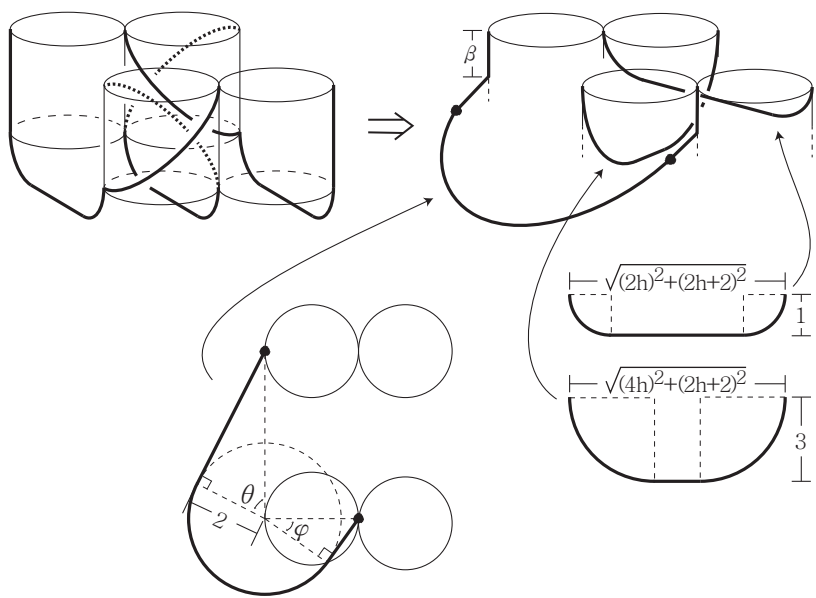}
\caption{Two crossings lying on other sides at the bottom area}
\label{fig10}
\end{figure}

Now consider the type that the crossing at the front lies on the left cylinder
and the crossing at the back lies on the right cylinder as shown in Figure \ref{fig10}.
We similarly have the same result for the type of its mirror reflection.
Delete all three subarcs of $K$,
and attach three new arcs as shown in the right figure.
One arc consists of two quarter circles of radius 1 and a horizontal line segment connecting them
with the total width $\sqrt{(2h)^2+(2h+2)^2}$,
that is the distance between two end-points to be connected.
Another arc consists of two quarter circles of radius 3 and a horizontal line segment connecting them
with the total width $\sqrt{(4h)^2+(2h+2)^2}$.
The third arc consists of a subarc of a circle with radius 2 lying on an $xy$-plane,
two horizontal line segments tangent to the circle at the two end-points of the circular arc,
and two vertical line segments of length some $\beta$ adjacent to the two horizontal line segments.
These vertical line segments are needed to maintain the distance between arcs at least 2.
$\beta= 0.1$ is enough for such purpose.
The center of the circle with radius 2 is away from the two end-points of the two vertical line segments
by $2h+2$ and $2h$.
The angles $\theta$ and $\varphi$ can be obtained from the equations
$\cos \theta = \frac{2}{2h+2}$ and $\cos \varphi = \frac{2}{2h}$.
The total length of these three arcs is
$(\sqrt{(2h)^2+(2h+2)^2} - 2 + \pi) + (\sqrt{(4h)^2+(2h+2)^2} - 6 + 3\pi)
+ (2(\frac{3}{2}\pi - \theta - \varphi) + \sqrt{(2h+2)^2 - 4} + \sqrt{(2h)^2 - 4} + 2\beta)$ ($\approx 27.19$).

\noindent {\bf Step 2.}
{\em Shortening the top part.\/}
\vspace{2mm}

\begin{figure}[h]
\includegraphics{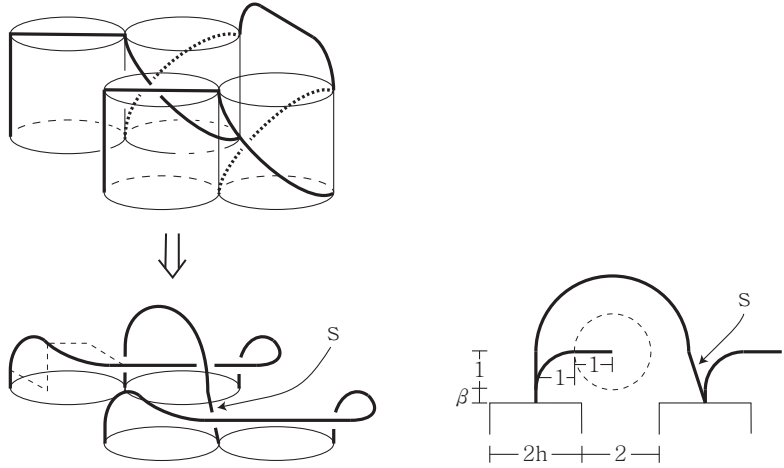}
\caption{$c(K)$ is even}
\label{fig11}
\end{figure}

First consider the case that $c(K)$ is even, that is, the top part consists of four cylinders of the same height.
We modify the top part as illustrated in Figure \ref{fig11}.
Delete all three subarcs of $K$ and attach three new arcs.
Each of two new arcs consists of a horizontal line segment of length $4h -2$,
four quarter circles of radius 1, and two vertical line segments of length $\beta = 0.1$.
More precisely, the horizontal line segment is parallel to the $y$-axis,
and each pair of quarter circles are connected and
one quarter circle of the pair lies on an $xy$-plane and the other lies on an $xz$-plane.
The third arc consists of a vertical line segment of length $1 + \beta$,
a half circle of radius 2, and  finally an almost vertical line segment $s$.
The total length of these three arcs is
$2(2\beta + 2\pi + 4h - 2) + (1 + \beta + 2\pi + \sqrt{(2h-2)^2+(1+\beta)^2})$ ($\approx 27.19$).

\begin{figure}[h]
\includegraphics{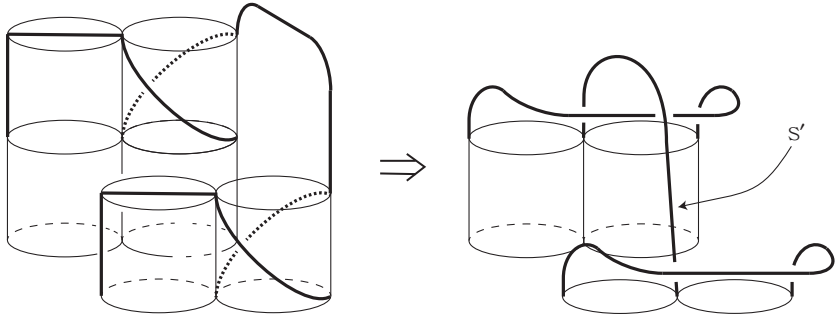}
\caption{$c(K)$ is odd}
\label{fig12}
\end{figure}

Now consider the case that $c(K)$ is odd, that is, the top two cylinders at the back are $2h$ higher than
the top two cylinders at the front as shown in Figure \ref{fig12}.
The construction is very similar to the case of even $c(K)$.
Only different thing is that the line segment $s$ is replaced by a longer line segment $s'$.
The total length of these three arcs is
$2(2\beta + 2\pi + 4h - 2) + (1 + \beta + 2\pi + \sqrt{(2h-2)^2+(1+\beta + 2h)^2})$ ($\approx 29.52$).

By using the mathematical software Mathematica,
we confirmed that each construction can be realized as a tube with uniform radius 1.
As a conclusion, the total ropelength of this embedding is bounded above by
$(2h\sqrt{\pi^2+4} + 2 \pi h+3 \pi + 2h) +
(6\pi + 8h + 5\beta - 3 + \sqrt{(2h-2)^2+(1+\beta + 2h)^2}) +
2h (\sqrt{\pi^2 + 4}+1)(c(K)-4)
\leq 11.39 c(K) + 12.37$.

This completes the proof of Theorem \ref{thm:Rop2}.


\begin{thebibliography}{ABGWQQ}
\bibitem {ACCJSZ} C. Adams, M. Chu, T. Crawford, S. Jensen, K. Siegel, and L. Zhang,
    {\em Stick index of knots and links in the cubic lattice},
    J. Knot Theory Ramifications \textbf{21} (2012) 1250041.
\bibitem {B} G. Buck,
    {\em Four-thirds power law for knots and links},
    Nature \textbf{392} (1998) 238--239.
\bibitem {BO} G. Buck and J. Orloff,
    {\em A simple energy function for knots},
    Topology Appl. \textbf{61} (1995) 205--214.
\bibitem {BS} G. Buck and J. Simon,
    {\em Thickness and crossing number of knots},
    Topology Appl. \textbf{91(3)} (1999) 245--257.
\bibitem {BZ} G. Burde and H. Zieschang,
    {\em Knots},
    Walter de Gruyter \& Co. (1985).
\bibitem {CFM} J. Cantarella, X. Faber, and C. Mullikin,
    {\em Upper bounds for ropelength as a function of crossing number},
    Topology Appl. \textbf{135} (2004) 253--264.
\bibitem {CKS1} J. Cantarella, R. Kusner, and J. Sullivan,
    {\em Tight knot values deviate from linear relations},
    Nature \textbf{392} (1998) 237--238.
\bibitem {CKS2} J. Cantarella, R. Kusner, and J. Sullivan,
    {\em On the minimum ropelength of knots and links},
    Invent. Math. \textbf{150(2)} (2002) 257--286.
\bibitem {C} J. Conway,
    {\em An enumeration of knots and links, and some of their algebraic properties},
    Computational Problems in Abstract Algebra (1970) 329--358.
\bibitem {D1} Y. Diao,
    {\em Minimal knotted polygons on the cubic lattice},
    J. Knot Theory Ramifications \textbf{2} (1993) 413--425.
\bibitem {D2} Y. Diao,
    {\em The Lower Bounds of the Lengths of Thick Knots},
    J. Knot Theory Ramifications \textbf{12} (2003) 1--16.
\bibitem {DE} Y. Diao and C. Ernst,
    {\em The Complexity of Lattice Knots},
    Topology Appl. \textbf{90} (1998) 1--9.
\bibitem {DEPZ} Y. Diao, C. Ernst, A. Por, and U. Ziegler,
    {\em The ropelengths of knots are almost linear in terms of their crossing numbers\/},
    Preprint arXiv:0912.3282v1.
\bibitem {DEY} Y. Diao, C. Ernst, and X. Yu,
    {\em Hamiltonian knot projections and lengths of thick knots},
    Topology Appl. \textbf{136} (2004) 7--36.
\bibitem {DET} Y. Diao, C. Ernst, and M. Thistlethwaite,
    {\em The linear growth in the length of a family of thick knots},
    J. Knot Theory Ramifications \textbf{12} (2003) 709--715.
\bibitem {DEZ} Y. Diao, C. Ernst, and U. Ziegler,
    {\em The linearity of the ropelengths of Conway algebraic knots in terms of their crossing numbers},
    Kobe J. Math. \textbf{28} (2011) 1--19.
\bibitem {HKNO} K. Hong, H. Kim, S. No, and S. Oh,
    {\em Minimum lattice length and ropelength of knots},
    To appear in J. Knot Theory Ramifications.
\bibitem {HNO1} K. Hong, S. No, and S. Oh,
    {\em Upper bound on lattice stick number of knots},
    Math. Proc. Camb. Phil. Soc. \textbf{155} (2013) 173--179.
\bibitem {HNO2} K. Hong, S. No, and S. Oh,
    {\em Upper bounds on the minimum length of cubic lattice knots},
    J. Phys. A: Math. Theor. \textbf{46} (2013) 125001.
\bibitem {HNRAV} X. Hua, D. Nguyen, B. Raghavan, J. Arsuaga, and M. Vazquez,
    {\em Random state transitions of knots: a first step towards modeling unknotting by type {II} topoisomerases},
    Topology Appl. \textbf{154} (2007) 1381--1397.
\bibitem {HO1} Y. Huh and S. Oh,
    {\em Lattice stick numbers of small knots},
    J. Knot Theory Ramifications \textbf{14} (2005) 859--867.
\bibitem {HO2} Y. Huh and S. Oh,
    {\em Knots with small lattice stick numbers},
    J. Phys. A: Math. Theor. \textbf{43} (2010) 265002.
\bibitem {JaP} E. Janse van Rensburg and S. Promislow,
    {\em Minimal knots in the cubic lattice},
    J. Knot Theory Ramifications \textbf{4} (1995) 115--130.
\bibitem {JiP} G. T. Jin and W. K. Park,
    {\em Prime knots with arc index up to 11 and an upper bound of arc index for non-alternating knots},
    J. Knot Theory Ramifications \textbf{19} (2010) 1655--1672.
\bibitem {K} L. Kauffman,
    {\em State models and the Jones polynomial},
    Topology \textbf{26} (1987) 395--407.
\bibitem {LSDR} R. Litherland, J. Simon, O. Durumeric, and E. Rawdon,
    {\em Thickness knots},
    Topology Appl. \textbf{91} (1999) 233--244.
\bibitem {Mc} L. McCabe,
    {\em An upper bound on edge numbers of 2-bridge knots and links},
    J. Knot Theory Ramifications \textbf{7} (1998) 797--805.
\bibitem {Mu} K. Murasugi,
    {\em Jones polynomials and classical conjectures in knot theory},
    Topology \textbf{26} (1987) 187--194.
\bibitem {SIADSV} R. Scharein, K. Ishihara, J. Arsuaga, Y. Diao, K. Shimokawa, and M. Vazquez,
    {\em Bounds for the minimum step number of knots in the simple cubic lattice},
    J. Phys. A: Math. Theor. \textbf{42} (2009) 475006.
\bibitem {T} M. Thistlethwaite,
    {\em A spanning tree expansion of the Jones polynomial},
    Topology \textbf{26} (1987) 297--309.

\end{thebibliography}
\end{document}